\documentclass[12pt,a4paper]{article}\usepackage[top=20mm,left=12mm]{geometry}\textwidth18.2cm\textheight25cm
\usepackage{amsmath,amssymb}
\usepackage{float}\floatstyle{plaintop}\restylefloat{table}
\usepackage[tableposition=top]{caption}
\usepackage{graphicx,pstricks,fancyvrb,multirow,lineno,url}
\usepackage{listings}\usepackage{paralist, longtable}
\newtheorem{theorem}{Theorem}[section]
\newtheorem{remark}[theorem]{Remark}
\def\bexe{\begin{exercise}}\def\eexe{\eex\end{exercise}}
\def\bsol{\begin{solution}}\def\esol{\eex\end{solution}}
\def\bexa{\begin{example}}\def\eexa{\end{example}}
\def\brem{\begin{remark}}\def\erem{\end{remark}}
\def\bthm{\begin{theorem}}\def\ethm{\end{theorem}}
\def\blem{\begin{lemma}}\def\elem{\end{lemma}}
\def\bcor{\begin{corollary}}\def\ecor{\end{corollary}}
\def\bdefi{\begin{definition}}\def\edefi{\end{definition}}
\newcommand{\IDEA}{\textbf{Idea of the Proof.} }

\def\bmip{\begin{minipage}{\textwidth}}\def\emip{\end{minipage}}
\def\huga#1{\begin{gather} #1 \end{gather}}

\def\hual#1{\begin{align} #1 \end{align}}

\newcommand{\R}{{\mathbb R}}

\def\setm{\setminus}

\def\ga{\gamma}

\def\pa{{\partial}}\def\lam{\lambda}
\newcommand{\bi}{\begin{itemize}}\newcommand{\ei}{\end{itemize}}
\newcommand{\ben}{\begin{enumerate}}\newcommand{\een}{\end{enumerate}}
\newcommand{\bce}{\begin{center}}\newcommand{\ece}{\end{center}}
\newcommand{\bci}{\begin{compactitem}}\newcommand{\eci}{\end{compactitem}}
\newcommand{\bcen}{\begin{compactenum}}\newcommand{\ecen}{\end{compactenum}}
\newcommand{\bcena}{\begin{compactenum}[(a)]}
\newcommand{\reff}[1]{(\ref{#1})}

\newcommand{\hs}[1]{{\hspace{#1}}}\newcommand{\vs}[1]{{\vspace{#1}}}

\newcommand{\barr}{\begin{array}}\newcommand{\earr}{\end{array}}
\newcommand{\bpm}{\begin{pmatrix}}\newcommand{\epm}{\end{pmatrix}}
\newcommand{\bsm}{\left(\begin{smallmatrix}}
\newcommand{\esm}{\end{smallmatrix}\right)}
\newcommand{\ba}{\begin{array}}\newcommand{\ea}{\end{array}}

\def\Om{\Omega}

\def\eex{\hfill\mbox{$\rfloor$}}

\def\Ga{\Gamma}

\def\bd{\begin{displaymath}} \def\ed{\end{displaymath}}
\def\ba{\begin{array}} \def\ea{\end{array}}

\def\pdep{{\tt pde2path}}

\def\mlab{{\tt Matlab}}

\def\dhome{/hh/path/pde2path/demos/acsuite}

\def\Gbc{G_{\text{BC}}}

\newlength{\tew}

\def\ig{\includegraphics}

\renewcommand{\arraystretch}{1.05}\renewcommand{\baselinestretch}{1.0}
\def\medskip{}\def\bigskip{}
\usepackage{setspace}
\usepackage{hyperref}
\usepackage[ho,hu,12]{optional}
\def\taskip{\renewcommand{\arraystretch}{1}\renewcommand{\baselinestretch}{1}}
\def\teskip{\renewcommand{\arraystretch}{1.1}\renewcommand{\baselinestretch}{1.1}}
\def\hulst#1#2{\taskip\lstinputlisting[#1]{#2}\teskip}
\def\hutab#1{\taskip\begin{\table}#1\end{table}\teskip}
\def\trulle{{\tt trullekrul}}\def\ma{mesh adaptation} 
\def\ina#1{}\def\ins#1{}\newfloat{Algorithm}{ht}{lop}[section]
\def\slfig{\opt{hu,ho}{\begin{figure}[ht]}\opt{sl}{\begin{figure}[H]}}
\def\slfigH{\opt{hu,ho}{\begin{figure}[H]}\opt{sl}{\begin{figure}[H]}}
\def\sltab{\opt{hu,ho}{\begin{table}[ht]}\opt{sl}{\begin{table}[H]}}
\def\dhome{trullem}

\input{lst}
\def\wtH{\widetilde{H}}
\begin{document}
\text{}\vspace{10mm} 
\bce\Large
Using {\tt trullekrul} in \pdep\ -- anisotropic mesh  adaptation for some 
Allen--Cahn models in 2D and 3D \\[4mm]
\normalsize 
Hannes Uecker\\[2mm]
\footnotesize
Institut f\"ur Mathematik, Universit\"at Oldenburg, D26111 Oldenburg, 
hannes.uecker@uni-oldenburg.de\\[3mm]
\normalsize
\today
\ece 
\begin{abstract} We describe by means of some examples 
how some functionality of the \ma\ package \trulle\ can be used in \pdep. 
\end{abstract}
\tableofcontents

\section{Introduction} \label{isec}
\def\Lelem{L_{\text{elem}}}\def\Llow{L_{\text{low}}}
\def\Lup{L_{\text{up}}}\def\Lmax{L_{\text{max}}}
The \mlab\ package \pdep\ \cite{p2pure, p2phome} is designed for numerical continuation and 
bifurcation analysis of systems of PDEs of the form 
\huga{\label{gformt}
M\pa_t u=\nabla\cdot(c\otimes\nabla u)-a u+b\otimes\nabla u+f, 
}
where $u=u(x,t)\in\R^N$, $t\ge 0$, $x\in\Om$, $\lam\in\R^p$ is a 
parameter vector, $M\in\R^{N\times N}$ is a mass matrix, which 
may be singular, and the coefficients 
$c,a,b$ and $f$ may depend on $x,\lam$ and $u$. 
Details on the terms in \reff{gformt}, the discretization of \reff{gformt} 
by the FEM, the boundary conditions associated to \reff{gformt},  
and how to use \pdep\ to compute branches of steady and 
time periodic solutions of \reff{gformt}, can be found in \cite{actut,pftut,hotheo}  
and some further tutorials which together with the software and demos 
can be downloaded at \cite{p2phome}.  

Here we explain how to use the anisotropic \ma\ package 
\trulle\ \cite{KEJ16, KEJ17} in \pdep. For this we 
extend the introductory tutorial \cite{actut} by a number of some advanced 
examples for steady Allen--Cahn type problems, i.e., problems of type 
\huga{\label{ac00} 
G(u):=-c\Delta u+f(u),\quad u=u(x)\in\R,\quad x\in\Om\subset\R^d 
\text{ a bounded domain}, 
}
with diffusion constant $c$, 'nonlinearity' (everything except diffusion) 
$f:\R\to\R$ (which also depends on parameters). 
For $\Om$ we shall restrict to rectangles 
(2D, $d=2$) and cuboids (3D, $d=3$), and $f$ and 
the boundary conditions will be given for the specific examples below. 

In \cite{actut} we also explain some settings for \ma\ for problems 
of type \reff{ac00} based on mesh refinement 
using standard a posteriori error estimators, in 1D 
and 2D. In the context of solution branches $s\mapsto(u(s),\lam(s))$ for \reff{ac00} 
as computed by the main user--interface routine {\tt cont} of \pdep, 
where $\lam$ is used as a symbol for the active continuation parameter, 
it is desirable to adapt the mesh during continuation, for instance 
 after a number of continuation steps, 
or if the error estimate is above some user given threshhold. 
This mesh {\em adaptation} so far has 
been done in \pdep\ in a simple ad hoc way: first we {\em coarsen} the current 
mesh to a given (essentially fixed and uniform) background mesh  by interpolation of 
the current solution to the background mesh, and then generate a new 
mesh by {\em refining} the coarse mesh and solution. The coarsening 
is necessary because otherwise we may end up with (at places) unnessarily fine 
meshes. However, the simple approach sketched above is not very 
efficient, may lead to undesired branch--switching after coarsening, 
and, moreover, has so far not been fully implemented in 3D. 
By interfacing \pdep\ with \trulle, we now have genuine 2D and 3D {\em adaptation} 
options during continuation of branches, which means coarsening (only) where 
appropriate together with moving of mesh--points and refinement. 

In \S\ref{trsec} we explain 
the general setup of \trulle\ in \pdep, based on the standard data structures 
of \pdep\ with all data stored in the  struct {\tt p}, i.e., 
function handles to the rhs of \reff{ac0} and its Jacobian, 
FEM mesh, file-names for saving, controls for plotting, numerical constants 
such as stepsizes and solution tolerances), and the solution {\tt p.u} itself 
and the current tangent to the solution branch. 
For this background, the \pdep\ data structures, the continuation algorithms, 
and the general usage of \pdep, 
we refer to \cite{qsrc} and the tutorials at \cite{p2phome}, in particular 
\cite{actut}. In \S\ref{exsec} we then explain the usage of \trulle\ 
by some example problems of type \reff{ac0}. 

\section{General setup of \trulle\ in \pdep}\label{trsec}
 Given a function $u:\Om\to \R$, \trulle\ 
aims to optimize the FEM mesh to minimize the interpolation error 
$\|u-u_h\|_{L^p(\Om)}$, where $u_h$ as usual is the continuous piecewise linear 
interpolation of the nodal values. 
This is based on the discrete Hessian $H(u_h)$ of 
$u_h$, and the associated metric field 
\huga{\label{trpsi}
\Psi=\frac 1 \eta(\det(\wtH))^{-\frac 1 {2p+d}}\wtH, 
}
where $\wtH$ denotes a matrix of absolute values of the 
eigenvalues of $H(u_h)$, and $\eta$ is a scaling 
factor, which can be used to control the number of mesh points. 
In \cite{CSX07} it is shown that meshes that are uniform 
wrt to $\Psi$ minimize $\|u-u_h\|_{L^p(\Om)}$. Thus, 
using $\Psi$, edge lengths $\Lelem$ of the triangulation 
are computed wrt to the metric $\Psi$, 
and then elements are coarsened/refined according 
to the following algorithm (cf.~\cite[Fig.6]{LSB05}): 
\begin{Algorithm}
\ina{Algorithm!\trulle}
\fbox{\parbox{0.98\textwidth}{
Compute the maximum edge length $\Lmax$ in metric space, and perform 
the following loop until $\Lmax<\Lup$ or until {\tt it=imax}: 
\ben
\item Eliminate the edges shorter than $\Llow$ by coarsening the 
respective elements; this includes swapping of elements with wrong orientation. 
\item Refine the elements with $\Lelem>\Lup$ by splitting 
their longest edge and postprocessing (splitting of adjacent elements); 
\item Move mesh--points according to a smoothing algorithm based on the 
(discrete) Laplacian. 
\item Update $\Psi, \wtH$ and hence $\Lelem$ and $\Lmax$. 
\een 
}}
\caption{{\small Outline of \trulle\ mesh adaptation. Given a current solution $u_h$, 
the inner adaptation loop is performed, and the resulting new $u_h$ is 
used as an initial guess for a Newton loop on the new mesh. This outer 
step is repeated {\tt p.nc.ngen} times.} \label{trullealg}}
\end{Algorithm}

\teskip
The details of each of the steps 1 to 3 in Algorithm \ref{trullealg} 
are controlled by a number  of \trulle\ parameters, which in \pdep\ 
we store in the field {\tt p.trop} ({\tt tr}ullekrul {\tt op}tions) 
of the basic data struct {\tt p}.  
To initialize 
{\tt trop} we provide the two functions {\tt troptions2D} 
and {\tt troptions3D}, which as variants of the \trulle\ 
function {\tt gen\_options} 
set \trulle\ parameters in the way which appear to be most robust and 
efficient in 2D, respectively 3D. For us, the most important parameter 
is $\eta$ in \reff{trpsi}, where, for fixed $\Lup$, 
larger $\eta$ gives less triangles to refine. For adaptation during 
continuation we typically want to keep the number of mesh points $n_p$ 
below some bound, and thus $\eta$ should be set by some function which depends 
on $n_p$. To give maximum flexibility, the user may provide such a function in 
{\tt trop.etafu} with signature {\tt eta=etafu(p,np)}. 
The default setting is {\tt trop.etafu=@stanetafu}, which returns 
the constant $\eta=10^{-3}$, while a simple dependence on $n_p$ 
is given by {\tt eta=etafua(p,np)} which (by default) returns $10^{-5}{\tt np}$. 
Altogether, in Table \ref{trtab2} we list 
the {\tt trop} parameters/function handles which we find most useful for tuning 
the adaptation, and we strongly 
recommend to experiment with these and the other \trulle\ parameters, see 
the sources for detailed comments. 
The \pdep\ parameter to switch on 
mesh--adaptation by \trulle\ instead of option (i)  
within {\tt cont} is {\tt p.sw.trsw}$>0$. 

\begin{table}[ht]\taskip
\caption{{\small Most important parameters/function handles 
in {\tt p.trop}, 
see also 
{\tt troptions2D} and {\tt troptions3D} for further parameters and comments. 
{\tt p.trcop} (see below) 
needs the same parameters as {\tt p.trop}, and additionally 
{p.trcop.sw}, see bottom of table. Typically, we just copy {\tt trop} 
to {\tt trcop}, set {\tt p.trcop.npb}, and 
reset selected parameters such as {\tt p.trcop.sw}. 
\label{trtab2}}}
\bce\vs{-4mm}
{\small 
\begin{tabular}{p{0.12\tew}|p{0.87\tew}} 
Parameter&meaning, default value\\
\hline
etafu&default {\tt eta=stanetafu(p,np)} returning 0.001. 
Larger eta gives less elements to adapt. See also 
{\tt etafua}, yielding $\eta=10^{-5}n_p$, and we 
recommend to copy {\tt etafua} to the local directory and 
experiment with the prefactor.    \\
zfu&function handle with signature {\tt z=zfu(p)}, 
to select the field $u$ for which the interpolation error 
is estimated. Default {\tt zfu=@stanzfu} for which $z=u_1$ (first component 
of current solution). May be useful to overload for multi--component 
systems. In some cases, also scaling of $z$ is useful, e.g., 
$z=\exp(u_1)$. \\
setids&function handle 
needed to link the \pdep\ data structures with {\tt trullekrul} in case 
that different boundary segment numbers (and different BCs on different segments) are assigned in \pdep. Defaults: {\tt @setidssq} in 2D (corresponding to 
a rectangular domain, as also used in {\tt stanpdeo2D}); 
{\tt @setidsbar} in 3D (corresponding to 
a cuboid domain, as also used in {\tt stanpdeo3D}). \\
\hline 
sw&behavior of {\tt tradapt} according to Table \ref{trtab1}; default 15. \\
ppar&$p$ to optimize the interpol.~error in the $p$-norm metric. Default: 1000, i.e., close to $\infty$ norm.\\
innerit&number of iterations in {\tt trullekrul}, default 2. (Not to confuse with {\tt p.nc.ngen} giving the 'outer' number of adaptation iterations).\\
Llow, Lup&lower/upper thresholds for edges in metric space, defaults $1/\sqrt{2}$ and 
$\sqrt{2}$. Smaller 
Llow can be used to avoid too much coarsening.\\
qualP&weight for combining mesh qualtity in metric space and euclidean space. 
Defaults: 0 in 2D (metric space only), 2 in 3D (avoiding too acute tetrahedra) \\
\hline
trcop.npb&desired number of mesh-points for pure coarsening steps. \\
trcop.crmax&maximum number of pure coarsening steps.
\end{tabular}
}
\ece
\vs{-4mm}
\end{table}
\teskip
In some examples (in particular in 3D, see below), it is useful for 
mesh adaptation within {\tt cont} to first use only the coarsening and moving 
functionality of \trulle, and then adapt again with refinement. 
For this two step 
approach we provide a modification of the 
main \trulle\ wrapper function {\tt adapt\_mesh}, named {\tt tradapt}. 
Both,  {\tt adapt\_mesh} and {\tt tradapt}, 
take the option field {\tt trop} as last argument, and 
for {\tt tradapt} this should contain the field {\tt trop.sw} which 
encodes the operations 
\huga{\label{tr4}
\text{(face or edge) swapping, coarsening, moving (of mesh points), and refinement,} 
}
in a binary way  according to Table \ref{trtab1}. Thus, 
coarsening and moving as a preparatory step%
\footnote{executed if {\tt p.trcop.npb>0} and {\tt p.trcop.crmax>0} by calling 
{\tt tradapt(\ldots,p.trcop)}, i.e., with the '\trulle\ coarsening options' 
{\tt p.trcop} instead of the '\trulle\ options'  {\tt p.trop}}
 before adaptation 
is encoded as {\tt p.trcop.sw=5}, and additionally 
{\tt  p.trcop.npb} should contain the desired maximum number of 
meshpoints. 

\taskip 
\begin{table}[ht]\caption{{\small Control of {\tt tradapt} via {\tt sw} based binary 
coding with the 1st/2nd/3rd/4th bit switching on 
moving/refinement/coarsening/swapping, abbreviated as m/r/c/s, respectively. 
{\tt sw=15} thus corresponds to the original \trulle\ {\tt adapt\_mesh} 
behaviour. \label{trtab1}}}
\bce\vs{-4mm}
{\small 
\begin{tabular}{l|cccccccc|cccccccc} 
{\tt sw}&0&1&2&3&4&5&6&7&8&9&10&11&12&13&14&15\\
{\tt action}&none&m&r&m,r&c&c,m&c,r&c,r,m&s&m,s&r,s&m,r,s&c,s&c,m,s&c,r,s&c,r,m,s
\end{tabular}
}
\ece
\vs{0mm}
\end{table}
\teskip

\section{Example implementations and results}\label{exsec}
To illustrate the use and performance of the \trulle\ \ma\ we discuss 
some demos from {\tt pdepath/acsuite}, which also collects the demos 
discussed in \cite{actut}, and again we stress that new users should 
at least briefly browse \cite[\S4]{actut} and the associated demos. 

\subsection{2D} 
\paragraph{Extending {\tt ac2D} from \cite{actut}.} 
We start with \reff{ac00} on the rectangle 
$\Om=(-2\pi,2\pi)\times(-\pi,\pi)$, 
with Dirichlet BC, i.e. 
\begin{subequations}\label{acdbc}
\hual{\label{ac0}
&G(u):=-c\Delta u-\lam u-u^3+\ga u^5=0 \text{ in $\Om$}, \\
\label{2dbc} 
&\text{$u=d\cos(y/2)$ on $\Ga_2:=\{x=2\pi\}$, parameter $d$, and $u=0$ on 
$\pa\Om\setm\Ga_2$}. 
}
\end{subequations}
The label $\Ga_2=\{x=2\pi\}$ is due to the \pdep\ convention that 
the boundaries of rectangles as generated by {\tt stanpdeo2D} have the 
order bottom--right--top--left. 
For $d=0$, \reff{acdbc} features the bifurcation points 
\huga{
\lam_{jl}=(j/4)^2+(l/2)^2,\ \ \phi_{jl}=\sin(j(x+2\pi)/4)\sin(l(y+\pi)/2), \ \ 
j,l=1,2,\ldots 
} 
from the trivial branch $u\equiv 0$. The associated bifurcating branches 
have already been discussed in \cite[\S4.1]{actut} and are computed 
in {\tt ac2D/cmds1}. 
In {\tt ac2D/cmds2} we redo some of these computations 
starting with a very coarse mesh and aiming 
to illustrate the use and performance of \trulle. 
 Some results are shown in Fig.~\ref{f4b}, which compares the older \ma\ by 
error estimators with that by \trulle, for the solution on the 
first nontrivial branch at $\lam=4$ (see (a)). 
Clearly, {\tt e2rs} in (b) correctly identifies the triangles that are 
reasonable to refine, but refining the longest edges in Euclidean 
metric gives poor approximations at the 'boundary layers'. 
(c) shows a 'standard' refinement of (a) by {\tt adapt\_mesh} from 
\trulle\ (with the given parameters from {\tt troptions2D}), 
and (d) shows a refinement of (a) by {\tt tradapt} with {\tt sw=3}, 
i.e., without coarsening (and without swapping). 
This gives significantly more 
mesh points than (c), and our main purpose here is to illustrate the 
a--posteriori coarsening of (d) in (e,f), which gives a 
rather similar mesh as in (c). 

\slfig
\bce
\scalebox{0.9}{
\begin{tabular}{p{0.34\tew}p{0.34\tew}p{0.34\tew}}
{\small (a) Original mesh, $n_p{=}264$}&{\small (b) {\tt e2rs} refinement, 
$n_p{=}1200$}&
{\small (c) \trulle, $n_p{=}892$} \\
\ig[width=0.32\tew]{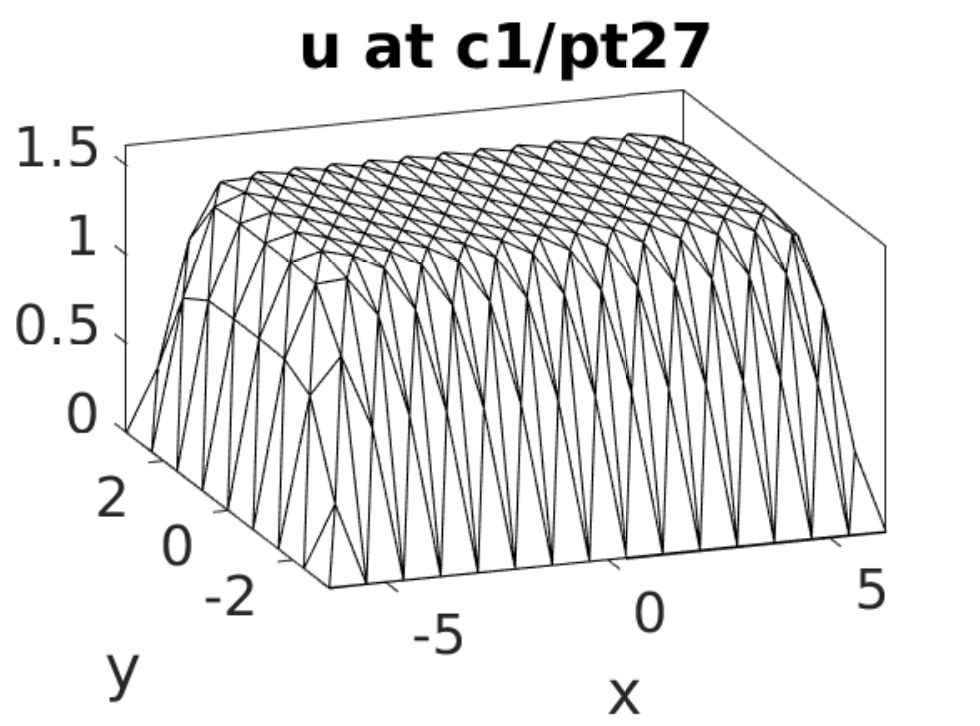}&\ig[width=0.32\tew]{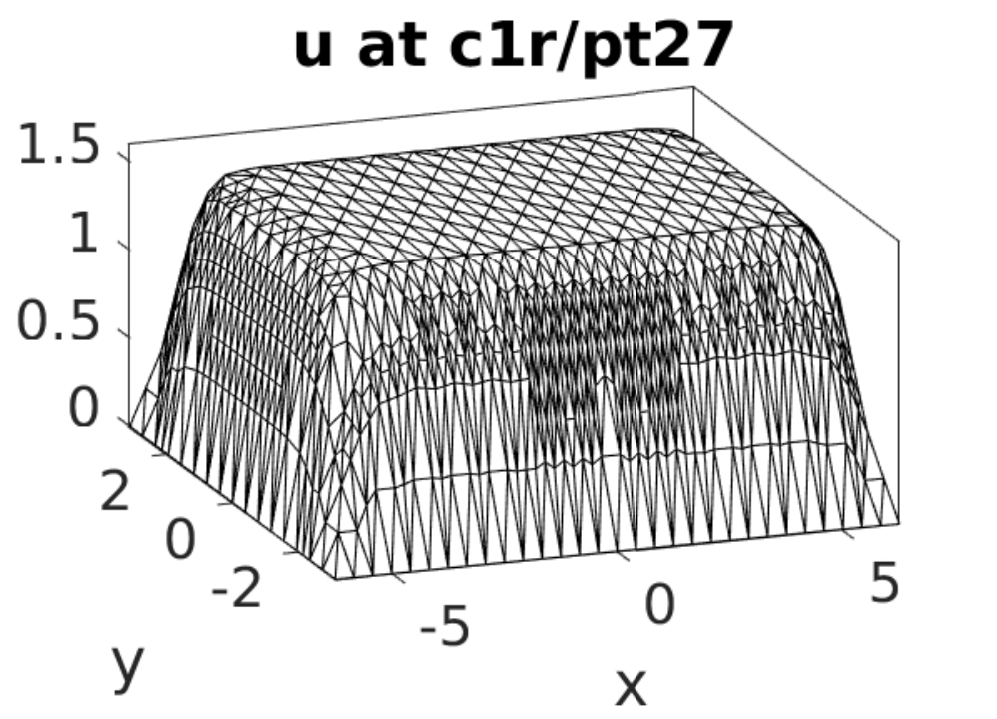}
&\ig[width=0.32\tew]{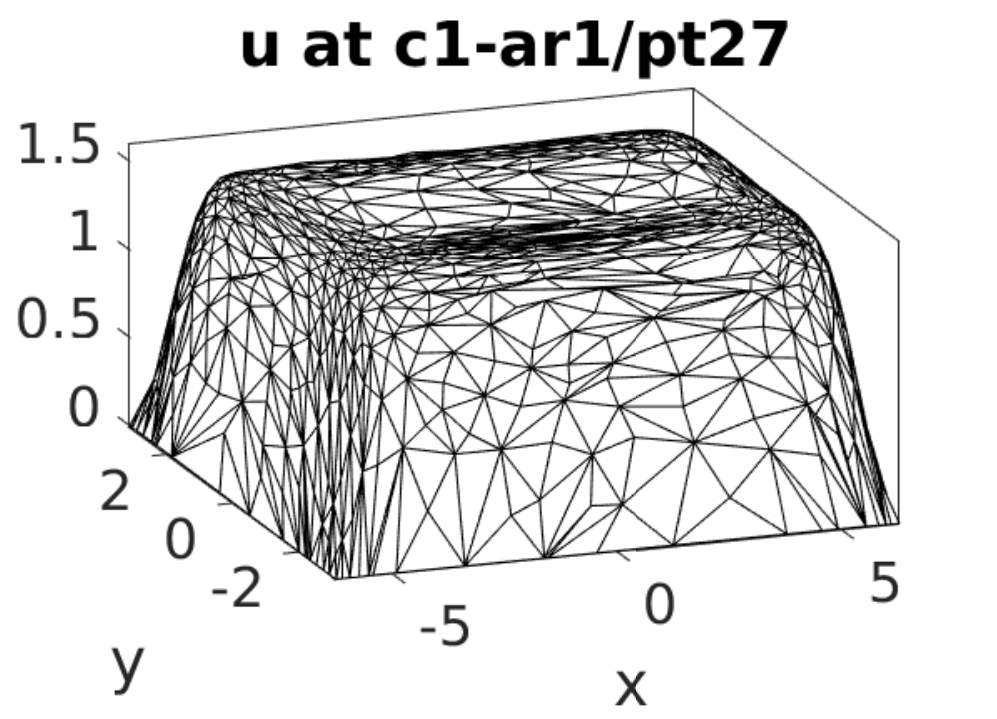}\\
\ig[width=0.32\tew]{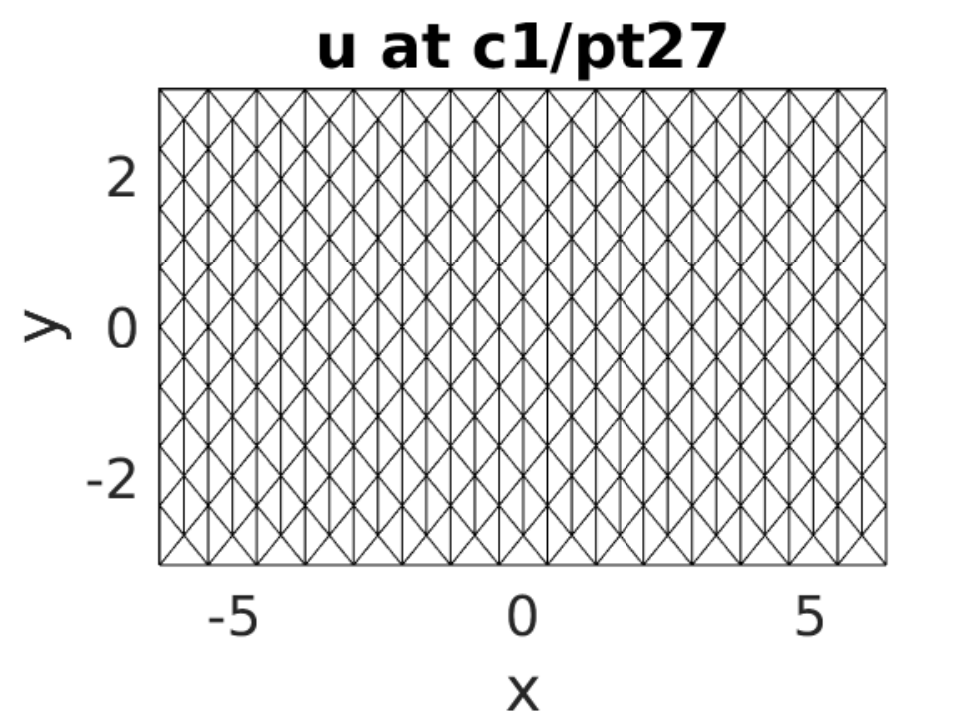}&\ig[width=0.32\tew]{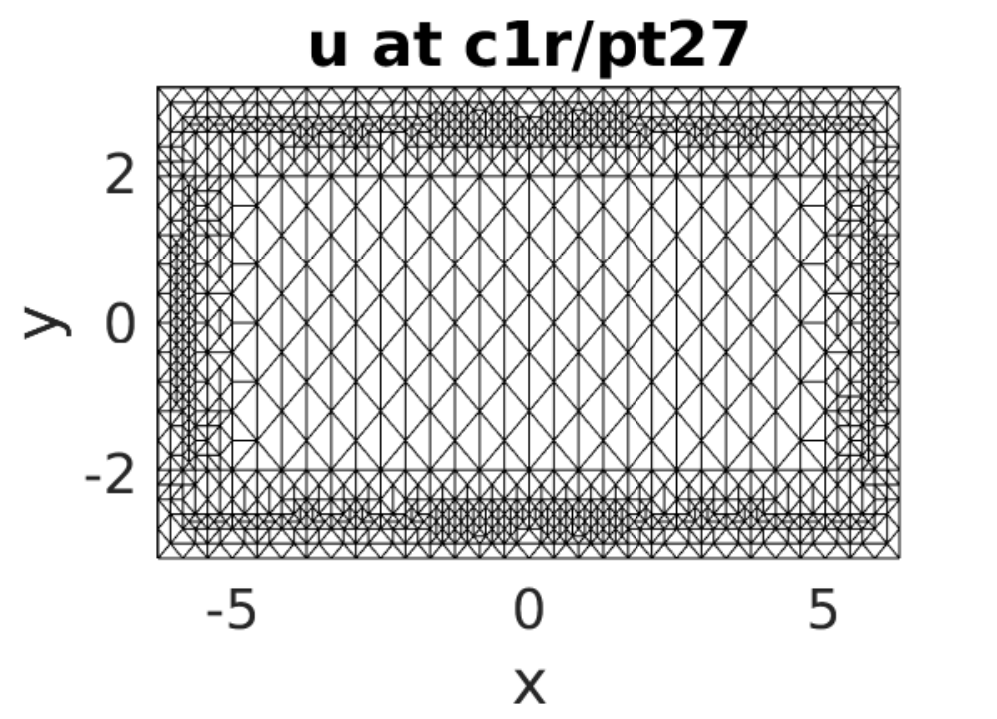}
&\ig[width=0.32\tew]{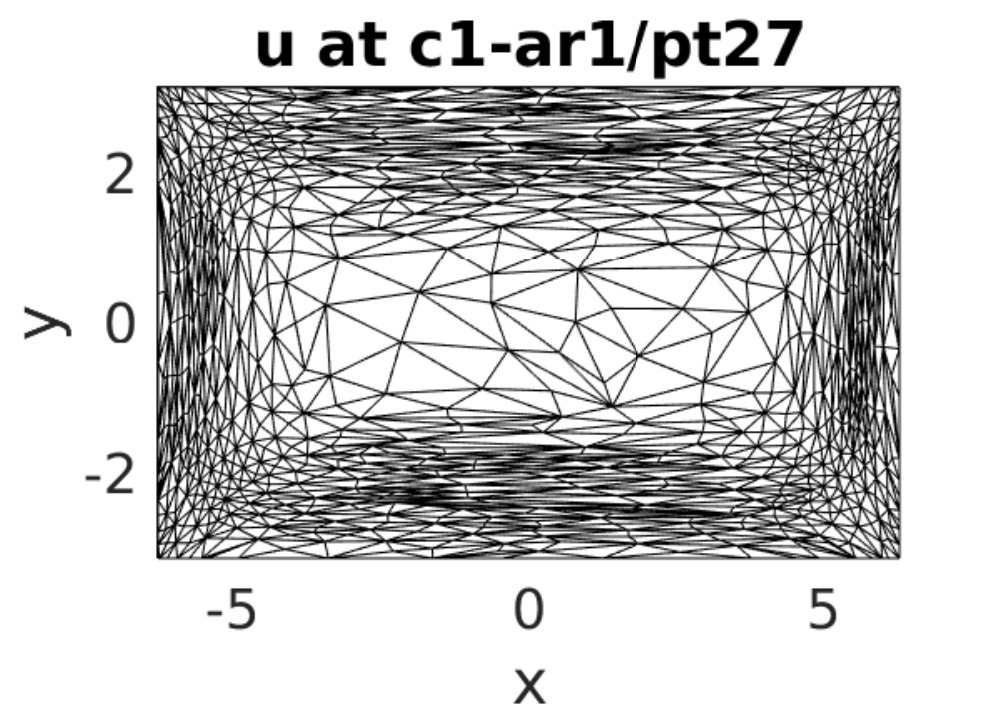}\\
{\small (d) no coarsening, $n_p{=}1616$}&
{\small \mbox{(e) a-posteriori coarsening $n_p{=}879$}}&{\small (f) top view of (d)}\\
\ig[width=0.32\tew]{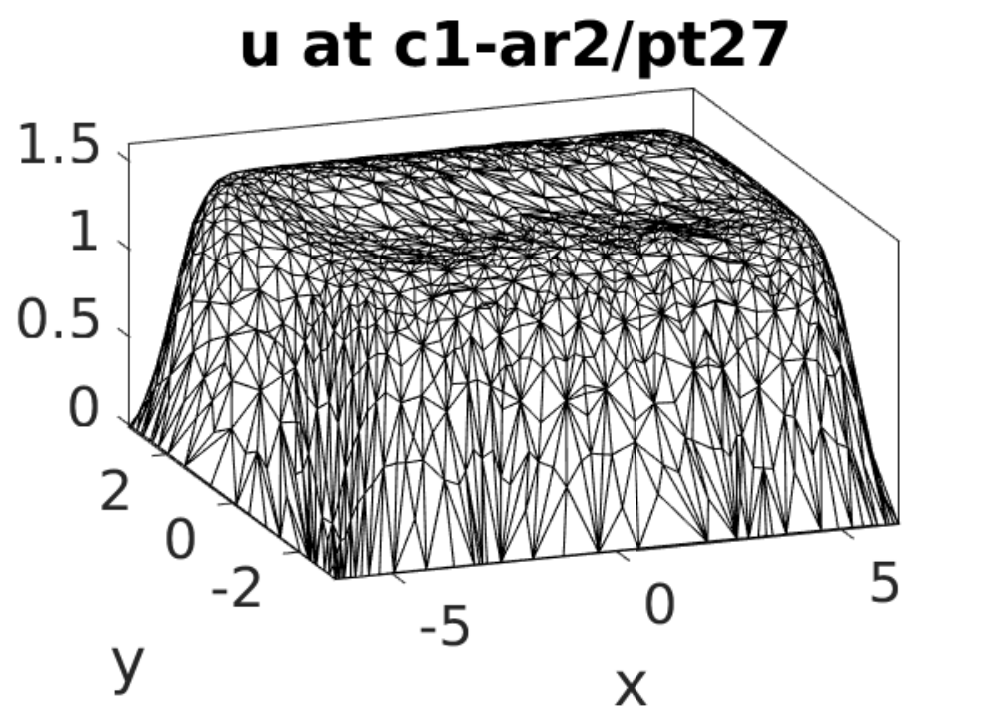}&\ig[width=0.32\tew]{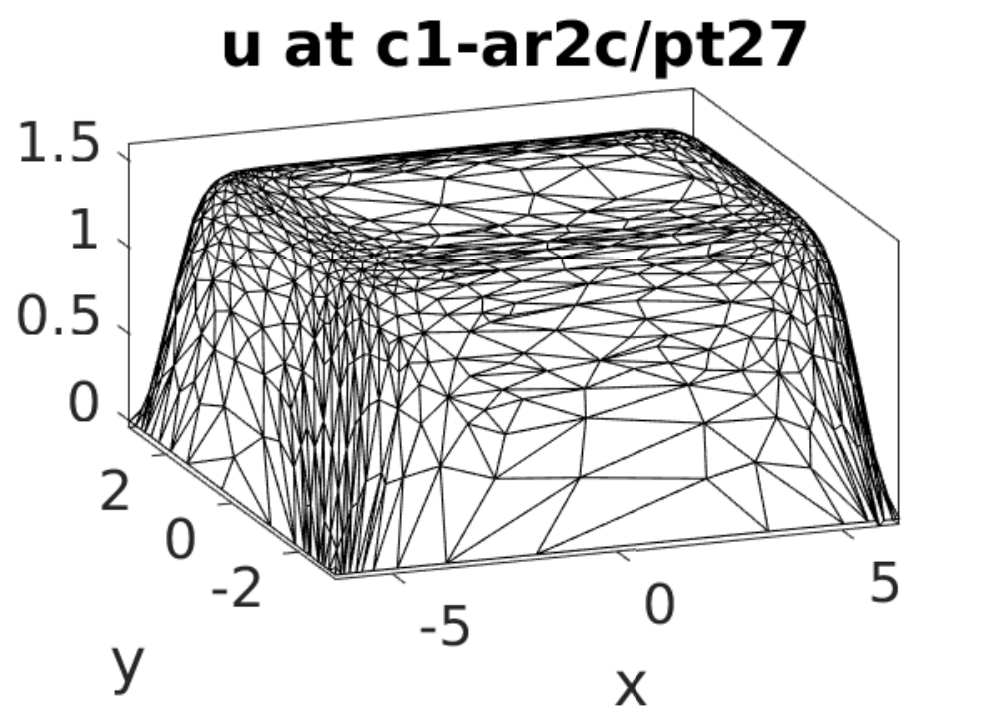}
&\ig[width=0.32\tew]{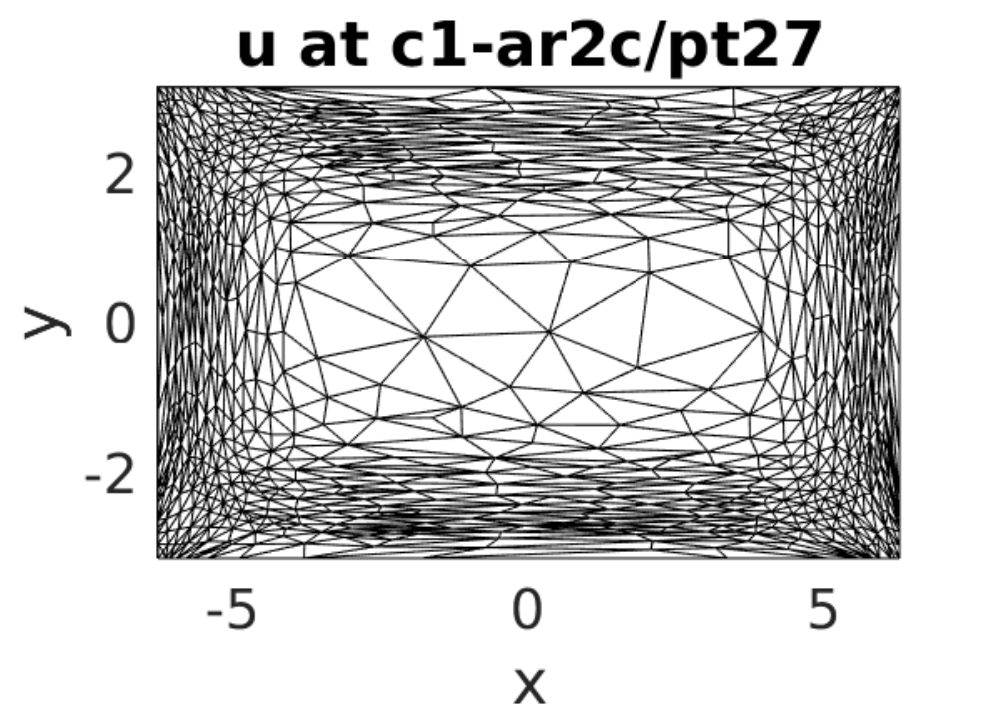}
\end{tabular}}
\ece 

\vs{-6mm}
   \caption{{\small Selection of plots generated in {\tt ac2D/cmds2.m}, 
illustrating different options for mesh-adaptation, see text for comments. 
  \label{f4b}}}
\end{figure}

\begin{figure}[ht] 
\bce
\scalebox{0.9}{
\begin{tabular}{p{0.28\tew}p{0.71\tew}}
\begin{tabular}{l}
{\small (a)}\\
\hs{-10mm}\ig[width=0.33\tew]{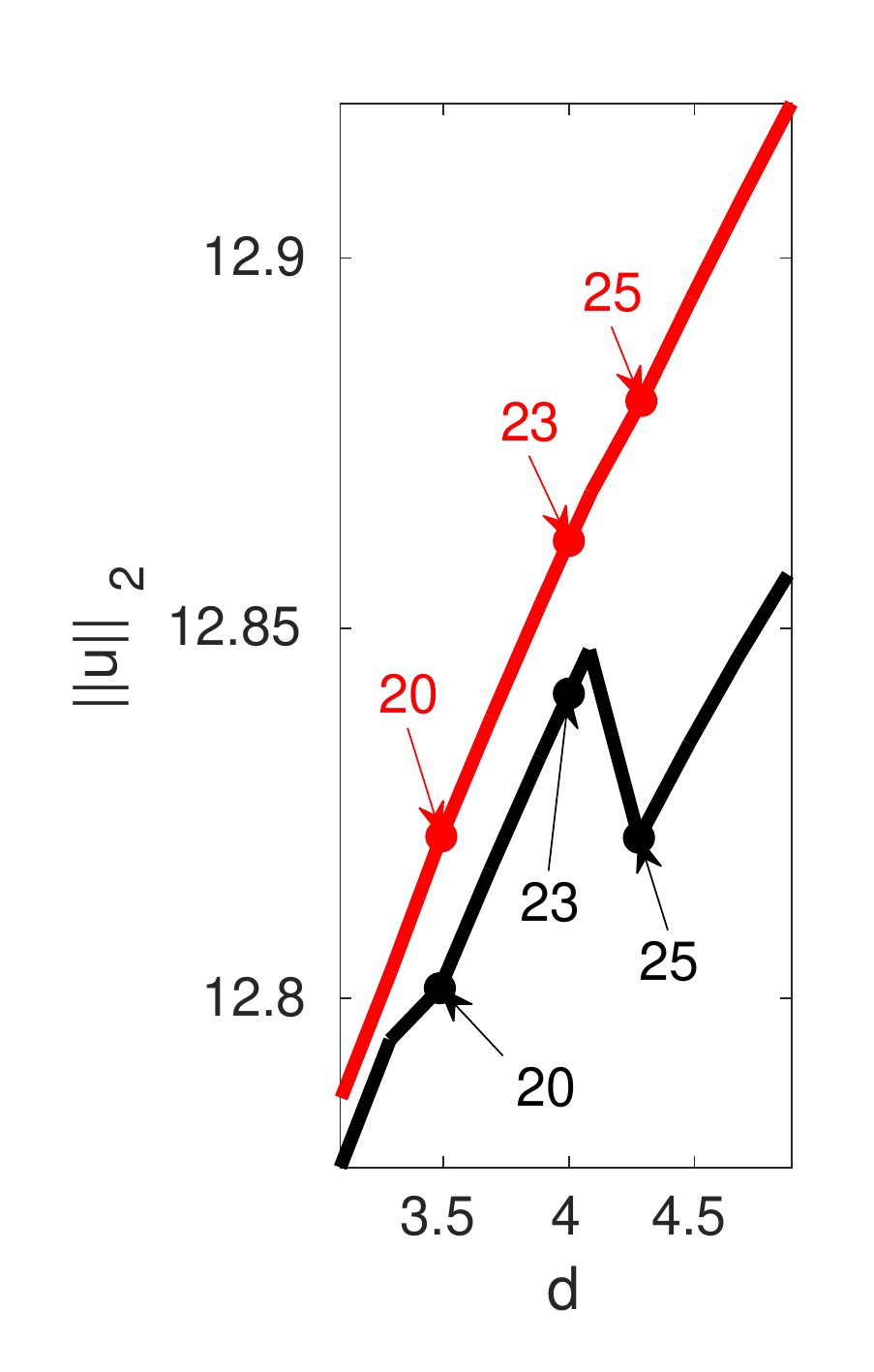}
\end{tabular}
&
\raisebox{-3mm}{\begin{tabular}{ll}
{\small (b) simple adaptation}\\
\hs{-3mm}\ig[width=0.4\tew]{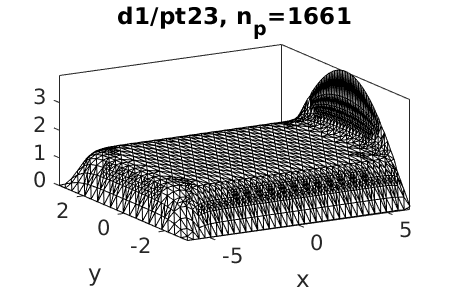}
&\hs{-3mm}\ig[width=0.4\tew]{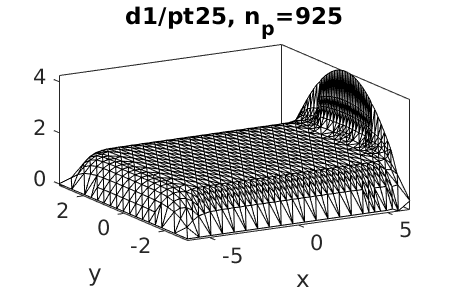}\\
{\small (c) \trulle\ adaptation}\\
\hs{-3mm}\ig[width=0.4\tew]{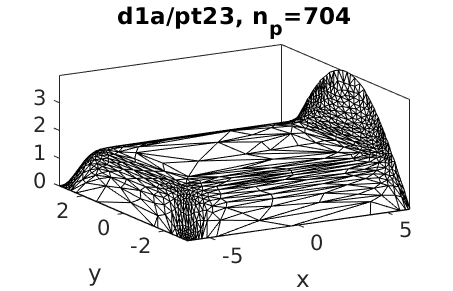}
&\hs{-3mm}\ig[width=0.4\tew]{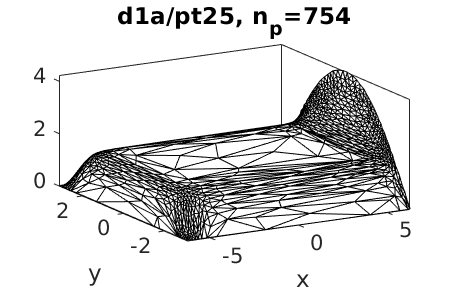}
\end{tabular}}
\end{tabular}
}
\ece 

\vs{-6mm}
   \caption{{\small Further results from 
{\tt cmds2}. Mesh adaptation every 5 steps of continuation in $d$, using 
{\tt e2rs} and \trulle, starting with the coarse mesh/solution 
from Fig.~\ref{f4b}(a). See text for further comments. 
  \label{f4c}}}
\end{figure}

Coarsening steps as in Fig.~\ref{f4b} from (d) to (e) 
are in particular useful for \ma\ during 
continuation. 
In Fig.~\ref{f4c} we continue the solution from Fig.~\ref{f4b}(a) (at $\lam=4$ 
and $d=0$) in $d$ with mesh adaption each 5th continuation step, 
again comparing \ma\ by \trulle\ with the old ad hoc adaptation by 
coarsening to the background mesh and then refining. 
The black branch in (a) belongs to the old option. 
The meshes and solutions (see (b) for two example plots) 
generally appear reasonable, but for larger $d$ the interpolation 
down to the coarse background mesh and subsequent refinement become 
problematic, i.e., both appear somewhat underresolved at the $x=2\pi$ 
boundary. A typical sign for such problems due to an inadequate background 
mesh are jumps in the branch data at adaption, that in (a) start to appear 
on the black branch for $d>3$. The controlled coarsening--refine approach 
by \trulle\ (red branch in (a)) is more robust in this respect. 
Of course, this is just one example, but it indicates a general result: 
if solutions during continuation develop boundary layers, or become 
strongly localized in some sense, use the \trulle\ mesh adaptation. 
Listing \ref{actl} shows some pertinent commands for this for easy review. 

\hulst{caption={{\small Commands from {\tt ac2D/cmds2.m} for 
continuation in $d$ with \ma\ by \trulle.}}, 
label=actl, language=matlab, stepnumber=5, firstnumber=42, linerange={42-53}}{\dhome/ac2D/cmds2.m}

\paragraph{A wandering boundary spot.} 
\def\dname{ac2Dwspot}
In a second example we consider \reff{ac00} on the 
rectangle $\Om=(-2\pi,2\pi)\times(-\pi,\pi)$ with the BCs 
\huga{\label{wsbc2}
\left\{\barr{l}
\text{$u=\exp(-(x-\xi)^2-z^2)$ on $\Ga_3=\{y=\pi\}$, parameter $\xi$,}\\
\text{$u=0$ on $\Ga_1=\{y=-\pi\}$, } \\
\text{$\pa_n u=0$ on $\pa\Om\setminus(\Ga_1\cup\Ga_3)$. }
\earr\right.
}
The purpose is to illustrate the mesh adaptation during continuation 
by considering a 'wandering spot' (upon continuation in $\xi$) on the top 
boundary. Additionally, 
this gives the opportunity to show how to put a parameter into the assembling 
of boundary values. 
The demo directory is {\tt \dname}, and Listing \ref{wsl1} shows the implementation 
of the rhs $G$, where we need to assemble the BCs in each step. 

\hulst{caption={{\small {\tt \dname/sGws.m}, $G$ for the wandering spot example. 
Due to the $\xi$ dependence, here we need to assemble the BCs $\Gbc$ in every 
call.  } }, 
label=wsl1, language=matlab,stepnumber=5}{\dhome/ac2Dwspot/sGws.m} 

\begin{figure}[ht]
\bce
\begin{tabular}{ll}
{\small (a)}&{\small (b)}\\
\ig[width=0.26\tew]{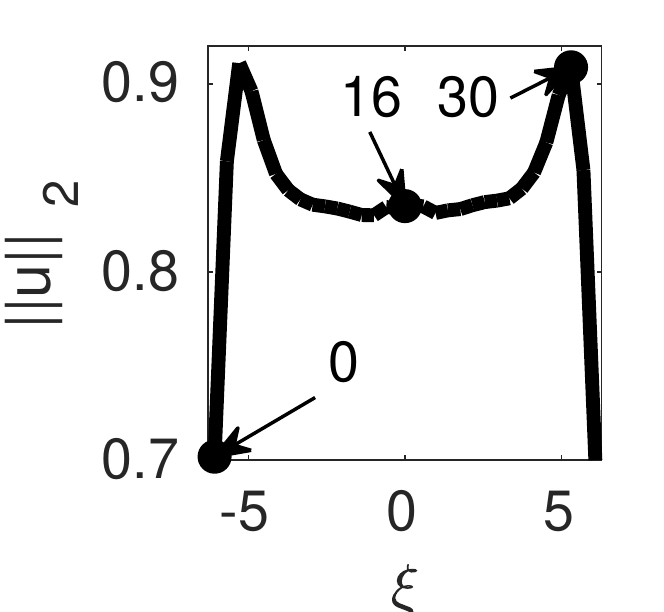}
&\ig[width=0.26\tew]{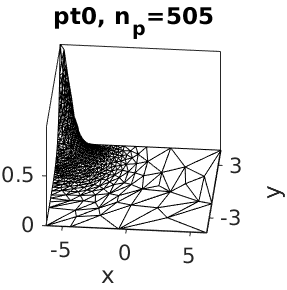}
\ig[width=0.26\tew]{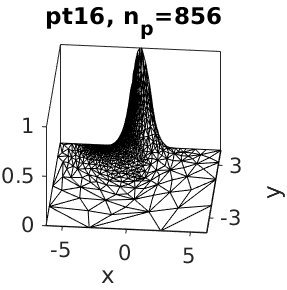}
\ig[width=0.26\tew]{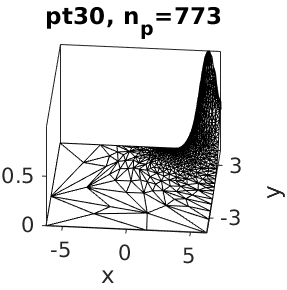}
\end{tabular}
\ece 
\vs{-4mm}
   \caption{{\small Example results from {\tt cmds1.m} for \reff{ac00} with the BCs \ref{wsbc2}, 
$(c,\lam,\ga)=(0.5,-0.25,1)$, continuation in $\xi$, with mesh--adapation 
every 5th step. 
 \label{wsf1}}}
\end{figure}

Fig.~\ref{wsf1} (see also Listing \ref{wsl2}) shows a continuation in $\xi$ 
in the subcritical case ($\lam=-0.25$), where the solutions are essentially 
characterized by the position of the boundary spot. For illustration 
we aim at rather coarse meshes, set {\tt amod=5} 
(mesh-adaptation every 5th step), and allow for extra coarsening steps in 
\trulle. The BD in (a) shows that this yields a reasonably smooth 
branch, and the sample plots in (b) that (as expected) the finest mesh 
may slighty lag behind the spot position (see in particular {\tt pt16}), 
but otherwise the coarsening/refinement setup works very well. 

\hulst{caption={{\small Selection from {\tt \dname/cmds1.m}. Continuation 
in $\xi$, $\lam=-0.25$ (subcritical case). } }, 
label=wsl2, language=matlab,stepnumber=5,firstnumber=19,linerange=19-24}{\dhome/ac2Dwspot/cmds1.m} 


In {\tt cmds2.m} we then continue in $\lam$ at fixed $\xi=0$. The 
black branch in (a) initially corresponds to the 'trivial' branch 
on which solutions are as in Fig.~\ref{wsf1} (specifically {\tt pt16} at $\xi=0$). 
Continuation in $\lam$ then yields an imperfect bifurcation to the 
primary unimodal branch in (a). Moreover, there are two BPs on the 
black branch, connected by the red branch. Alltogether, the \trulle\ \ma\ with 
{\tt amod=10} yields robust results with still rather coarse meshes. 

\hulst{caption={{\small Selection from {\tt \dname/cmds2.m}, continuation 
in $\lam$, $\xi=0$. We allow somewhat finer meshes, with 
${\tt etafub}=10^{-6}{\tt np}$. The remainder of {\tt cmds2} computes the bifurcating 
branch and then plots.  } }, 
label=wsl2, language=matlab,stepnumber=5,firstnumber=1,linerange=1-5}{\dhome/ac2Dwspot/cmds2.m} 

\begin{figure}[ht]
\bce
\scalebox{1}{
\begin{tabular}{l}
{\small (a) $\xi=0$, continuation in $\lam$ and sample plots from black branch}\\
\ig[width=0.24\tew]{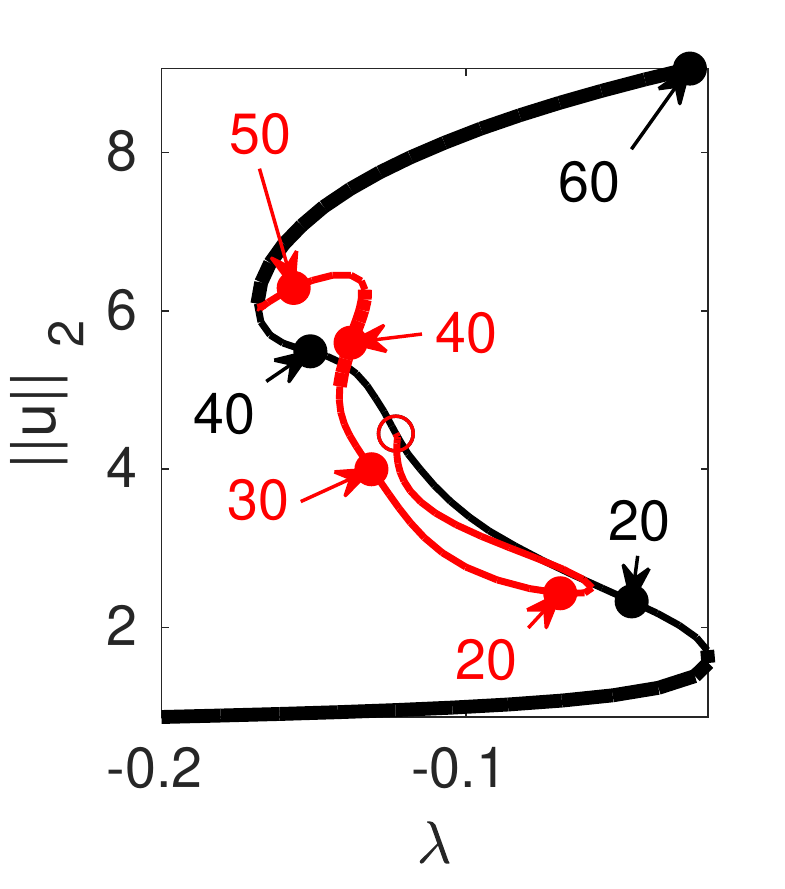}
\ig[width=0.25\tew,height=50mm]{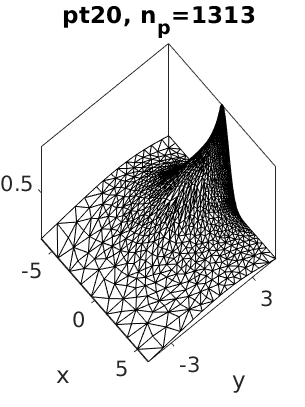}
\ig[width=0.25\tew,height=50mm]{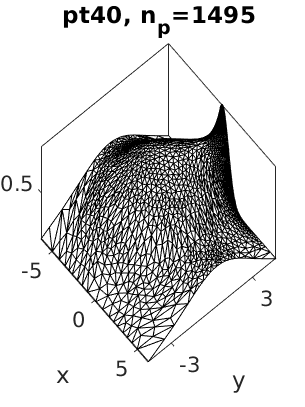}
\ig[width=0.25\tew,height=50mm]{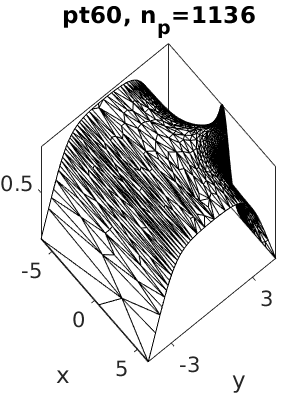}\\
{\small (b) Sample plots from red branch}\\
\ig[width=0.25\tew,height=50mm]{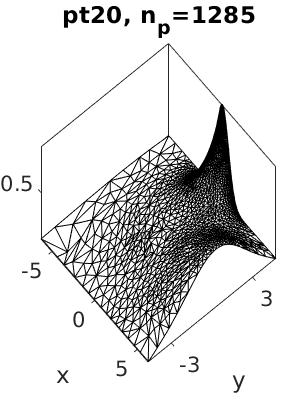}
\ig[width=0.25\tew,height=50mm]{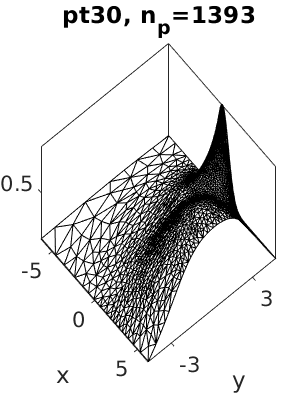}
\ig[width=0.25\tew,height=50mm]{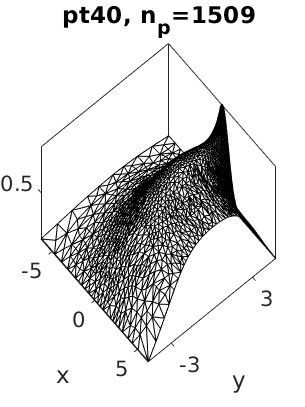}
\ig[width=0.25\tew,height=50mm]{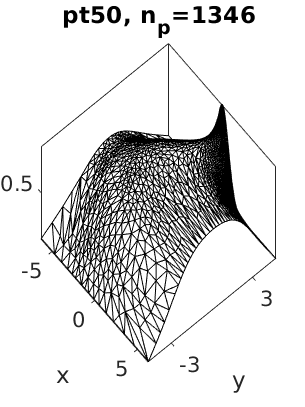}
\end{tabular}}
\ece 
\vs{-4mm}
   \caption{{\small Example from {\tt cmds2}. Continuing the primary solution branch 
in $\lam$, fixed $\xi=0$ yields an 'imperfect bifurcation' to the primary 
unimodal branch, and with an additional bifurcation to a branch with broken 
$x\mapsto -x$ symmetry, which, somewhat unexpectedly contains stable solutions, 
e.g., {\tt pt40}. 
 \label{acexws2}}}
\end{figure}

\subsection{3D}\label{ac3Dsec}
\begin{figure}[ht]
\bce
\begin{tabular}{p{0.33\tew}ll}
{\small (a) coarse initial mesh}&{\small (b) refinement of (a)}&
{\small (c) coarsening of (b)}\\
\hs{-6mm}\ig[width=0.35\tew]{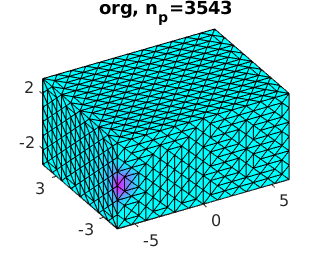}
&\hs{-4mm}\ig[width=0.35\tew]{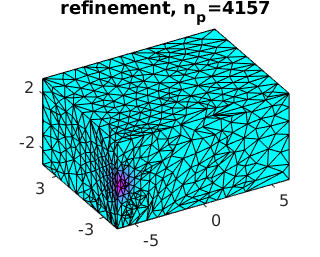}
&\hs{-4mm}\ig[width=0.35\tew]{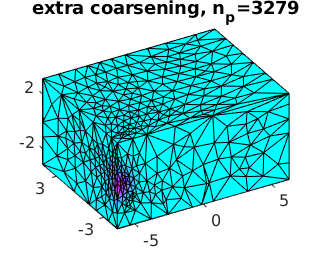}
\end{tabular}
\begin{tabular}{l}
{\small (d) BD, cont.~with fixed (red) and adaptive (blue) meshes, and 
sample solutions from the blue branch}\\
\hs{-6mm}\ig[width=0.33\tew,height=57mm]{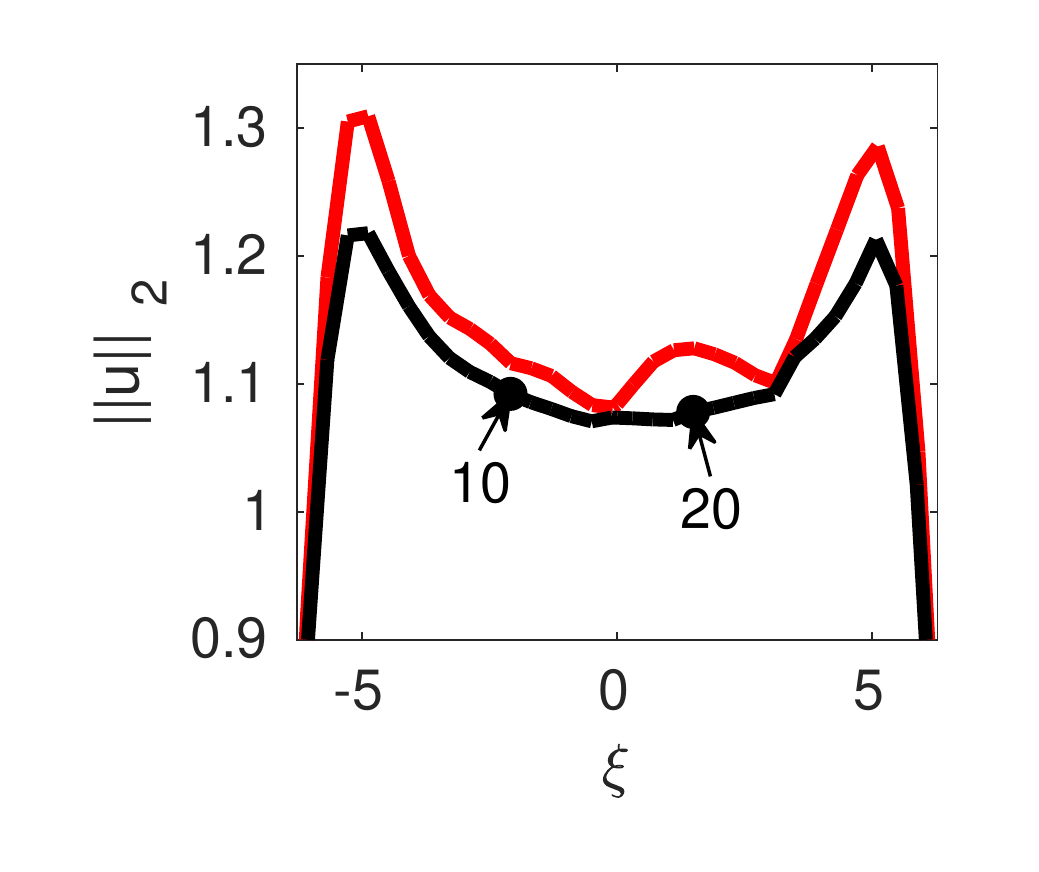}
\hs{2mm}\ig[width=0.35\tew]{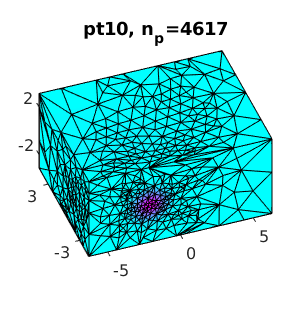}
\hs{2mm}\ig[width=0.35\tew]{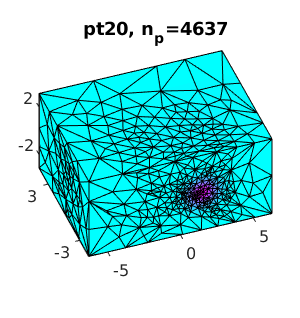}
\end{tabular}
\ece 
\vs{-8mm}
   \caption{{\small Selected plots from {\tt ac3Dwspot/cmds1}. Colorscale in all
 plots from $u=0$ to $u=1$.  \label{f5b}}}
\end{figure}

\def\dname{ac3Dwspot}
In {\tt \dname} we consider analogous BCs as in \reff{wsbc2}, i.e.%
\begin{subequations}\label{acdbc}
\hual{\label{ac0}
&G(u):=-c\Delta u-\lam u-u^3+\ga u^5=0 \text{ in $\Om$}, \\
\label{wsbc3}
&\left\{\barr{l}
\text{$u=\exp(-(x-\xi)^2-z^2)$ on $\Ga_3=\{y=-3\pi/2\}$, parameter $\xi$,}\\
\text{$u=0$ on $\Ga_5=\{y=3\pi/2\}$, } \\
\text{$\pa_n u=0$ on $\pa\Om\setminus(\Ga_3\cup\Ga_5)$.}
\earr\right.
}
\end{subequations}
For 3D cuboids as generated by {\tt stanpdeo3D}, the order 
of the boundary faces is bottom, left, front, right, back, top, which leads 
to a straightforward modification of {\tt sGws} from 2D to 3D. 
In {\tt cmds1} we again we start with a continuation in $\xi$ in the (here weakly) 
subcritical regime $\lam=0$, with the main results given in Fig.~\ref{f5b}, 
using a very coarse initial mesh 
with $n_p=3543$ to illustrate some effects and important settings 
for mesh--adaptation. 
Solutions on this coarse mesh look quite 
reasonable, see (a), and adaptive refinement (b) and 
subsequent coarsening (with {\tt p.trcop.sw=5}, see Table 
\ref{trtab1}) yield solutions which are visually the same, at least 
in the surface plots. However, 
for continuation on the fixed original coarse mesh from (a), 
for instance the $L^2$--norm  shows some 
unexpected fluctuations (red branch in (d)), e.g., for $\xi\in(0,2)$. 
These are essentially due to the spot 
at a given $\xi$ sitting on (close to) a mesh point, or in between 
two mesh points. The black branch in (d) is from continuation with 
mesh adaptation every 5th step, starting from the solution in (c), 
and setting {\tt p.trcop.npb=3000} to coarsen a given mesh to at 
most 3000 mesh points before refinement. The two sample plots on the 
right of (d) are from this blue branch. 
Some main observations are: 
\bcen
\item  The $L^2$--norm on the adapted meshes is generally slightly 
smaller than on the (coarse) uniform mesh, and 
\item The \trulle\ coarsening (to less than 3000 mesh points) --- refine approach 
yields nicely moving grids of small size with the finest meshing 
centered at the spots after each refinement. 
\item  Even though we only adapt every 5 steps, during which the spot 
moves to distance about 1 to 1.5 from the finest meshing, the 
black branch is reasonably smooth. However, if for instance we set {\tt amod=10} 
(such that the spot moves further away before remeshing), then visible jumps 
occur in the $L^2$--norm branch at each adaptation. 
\ecen 

As in 2D, in {\tt cmds2} we then switch to $\lam$ continuation, see Fig.~\ref{acexf2}. 
The basic behavior is as in 2D, i.e., the black branch turns into 
the primary unimodal (+spot) branch by an imperfect bifurcation, and 
there are now four bifurcation points connected by branches with broken 
$x\to -x$ symmetry.


\slfig
\bce
\scalebox{1}{
\begin{tabular}{l}{\small (a) BD and sample plots (lower left quarter 
domain) from the primary $\lam$ branch}\\
\ig[width=0.2\tew]{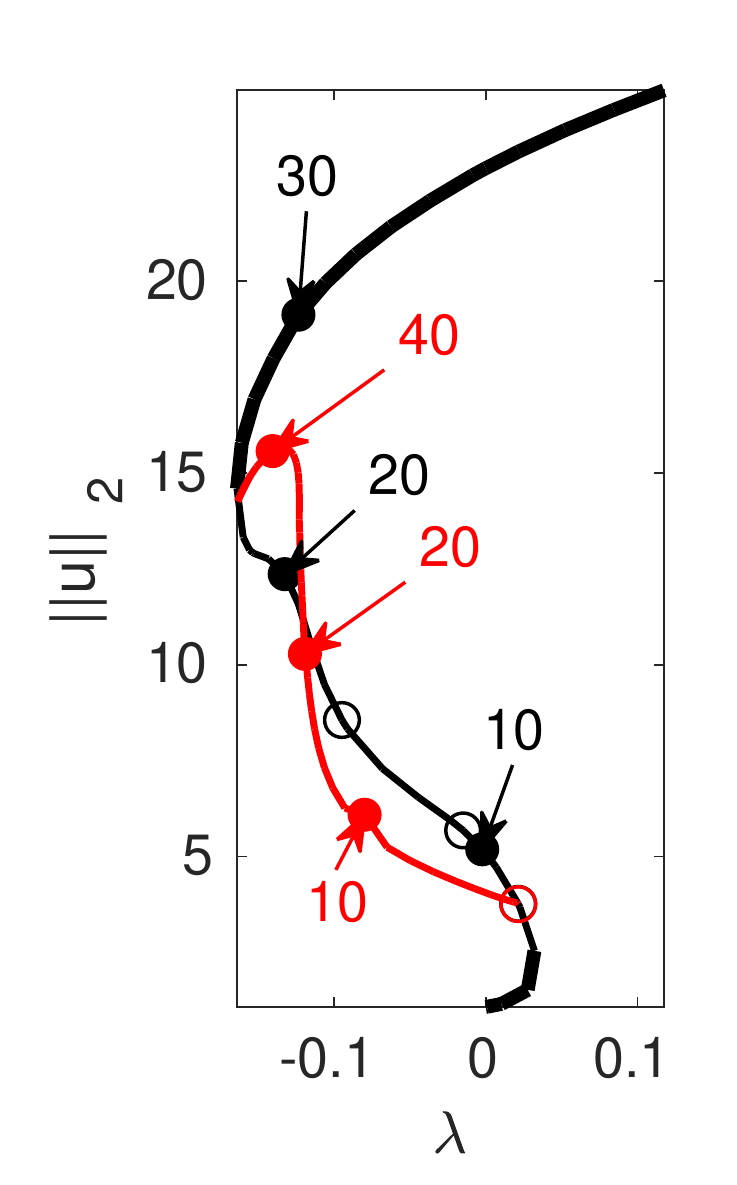}
\ig[width=0.29\tew]{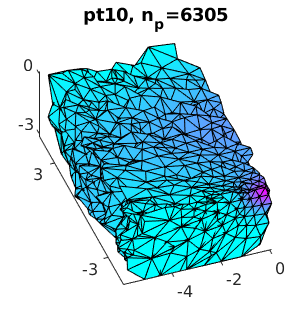}
\hs{-4mm}\ig[width=0.29\tew]{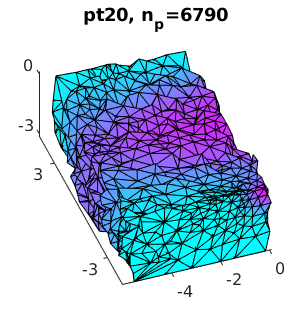}
\hs{-4mm}\ig[width=0.29\tew]{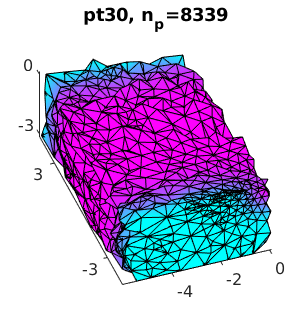}\\[3mm]
{\small (b) Sample plots (lower half of domain) from the first bifurcating branch (red)}\\
\hs{-10mm}\ig[width=0.41\tew]{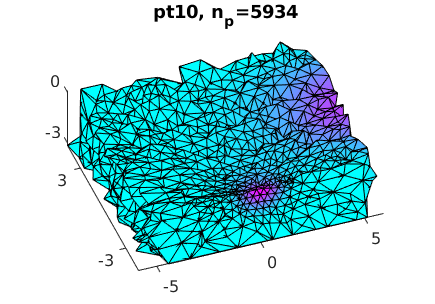}
\hs{-10mm}\ig[width=0.41\tew]{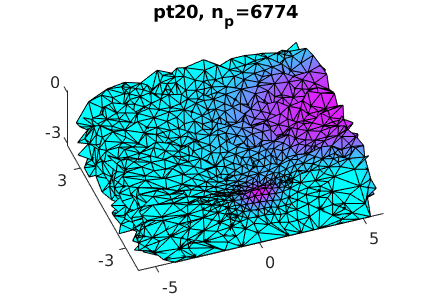}
\hs{-10mm}\ig[width=0.41\tew]{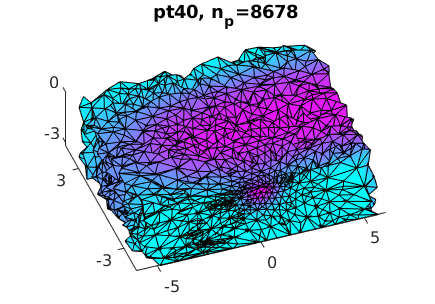}
\end{tabular}}
\ece 
\vs{-6mm}
   \caption{{\small Selected plots from {\tt ac3Dwspot/cmds2}. Colorscale in all
solution plots from $u=0$ to $u=1$. Fixed $\xi=0$, continuation in $\lam$, yielding imperfect 
bifurcation to primary unimodal branch (black).  At the first BP, a solution with 
broken $y\mapsto -y$ symmetry bifurcates, at the third BP 
the $z\mapsto -z$ symmetry is broken. 
The inner interfaces between $u=0$ and $u=1$ suggest a larger $n_p$. 
Some options used are: {\tt amod=8; trop.qualP=2.25; trcop.npb=8000}. 
 \label{acexf2}}}
\end{figure}


\renewcommand{\refname}{References}
\taskip
\small
\newcommand{\etalchar}[1]{$^{#1}$}

\end{document}